\documentclass[reprint,preprintnumbers]{revtex4-1} 

\usepackage{amsmath}  
\usepackage{amsfonts} 
\usepackage{graphicx} 
\usepackage{mathrsfs} 
\usepackage{bbm}
\usepackage{color}
\usepackage{hyperref}
\usepackage[cal=boondox,scr=boondoxo]{mathalfa}
\usepackage{mathtools}
\usepackage{bbold}

\begin{document}
\preprint{\textsl{KPOP}$\mathscr{E}$-2020-01}

\title{Proof of Cramer's rule with Dirac Delta Function
}

\author{June-Haak Ee}
\email{chodigi@gmail.com}
\affiliation{\textsl{KPOP}$\mathscr{E}$ Collaboration, Department of Physics, Korea University, Seoul 02841, Korea}

\author{Jungil Lee}
\email{jungil@korea.ac.kr} 
\thanks{Author to whom any correspondence should be addressed.
Director of Korea Pragmatist Organization for Physics Education (\textsl{KPOP}$\mathscr{E}$)}
\affiliation{\textsl{KPOP}$\mathscr{E}$ Collaboration, Department of Physics, Korea University, Seoul 02841, Korea}

\author{Chaehyun Yu}
\email{chyu@korea.ac.kr} 
\affiliation{\textsl{KPOP}$\mathscr{E}$ Collaboration, Department of Physics, Korea University, Seoul 02841, Korea}


\date{\today}

\begin{abstract}

We present a new proof of Cramer's rule by interpreting a system of linear equations 
as a transformation of $n$-dimensional Cartesian-coordinate vectors.
To find the solution, we carry out the inverse transformation by convolving the original coordinate vector with Dirac delta functions and 
changing integration variables from the original coordinates to new coordinates.
As a byproduct, we derive a generalized version of
Cramer's rule that applies to a partial set of variables, which is 
new to our best knowledge. 
Our formulation of finding a transformation rule for multi-variable functions
shall be particularly useful in changing a partial set of generalized coordinates
of a mechanical system. 

\end{abstract}

\maketitle 

\section{Introduction} 
Cramer's rule \cite{Cramer-1750} is a formula for solving a system of linear equations
as long as the system has a unique solution.
There are various proofs available \cite{Whitford-1953,Robinson-1970,Friedberg-2018,Brunetti-2014}
and six different proofs of Cramer's rule are listed in Ref.~\cite{Brunetti-2014}.

Dirac delta function $\delta(x)$ is not a proper function but a distribution
defined only through integrals~\cite{Dirac-1923,Dirac-1939,Dirac-1958}:
\begin{equation}
\label{del-def}
\int_{-\infty}^\infty dx\, \delta(x)f(x)=f(0),
\end{equation}
where $f(x)$ is a smooth function.
The Dirac delta function 
$\delta(x-a)$ projects out the value of a function $f(x)$ at a certain point 
$x=a$ after integration: $\int_{-\infty}^\infty dx\,f(x)\delta(x-a)=f(a)$. 
This elementary property of $\delta(x)$ can further be applied to change
the variable: $\int_{-\infty}^\infty dx\,f(x)\delta(x-y)=f(y)$. 
The Dirac delta function is particularly useful in physics.
One can make an educated guess that the change of variables can be made for a
multi-variable integral  by introducing the integration of products of
Dirac delta functions whose arguments contain the transformation constraints.

In this paper, we derive Cramer's rule by making use of the Dirac delta function to
change the variables of multi-dimensional integrals. 
The change of variables corresponds to a coordinate transform $\mathbbm{X}'=\mathbbm{R}\mathbbm{X}$, 
where $\mathbbm{X}=(x_1,\cdots,x_n)^T$ and $\mathbbm{X}'=(x_1',\cdots,x_n')^T$ are 
two $n$-dimensional Cartesian-coordinate column vectors and  $\mathbbm{R}=(R_{ij})$ is the corresponding invertible
transformation matrix. By changing the integration variables 
from $x_k'$'s to $x_k$'s in a trivial identity
$\mathbbm{X}=\prod_{k=1}^n\left[\int_{-\infty}^\infty dx_k'\delta(x_k'-\sum_{i_k=1}^nR_{ki_k}x_{i_k})\right]\mathbbm{X}$,
we carry out the inverse transformation $\mathbbm{X}=\mathbbm{R}^{-1}\mathbbm{X}'$.
To our best knowledge, this is a new proof of Cramer's rule.

Our formulation of finding a transformation rule for multi-variable functions
shall be particularly useful in changing a partial set of generalized coordinates
of a mechanical system and can further be extended to reorganizing the phase space of a multi-particle interaction
in particle phenomenology.

This paper is organized as follows. 
After providing a description of the notations that are frequently used in the remainder of this paper in Sec.~\ref{sec;def},
we proceed with the derivation of Cramer's rule by making use of Dirac delta functions in Sec.~\ref{sec;proof}. 
Section~\ref{sec;app} is devoted to demonstrate
how 
Cramer's rule for a partial set of variables
or coordinates
applies to a mechanical system.
In Sec.~\ref{sec;dis}, we discuss implications
and possible applications of our results.
A rigorous proof of Cramer's rule for a partial set of variables 
by employing mathematical induction
is presented in the Appendix.

\section{Definitions \label{sec;def}}
We define  $n$-dimensional column vectors 
$\mathbbm{X}$ and
$\mathbbm{X}'$
 whose $i$th elements are the Cartesian coordinates $x_i$ and $x_i'$, respectively:
\begin{equation}
\label{xxp}
\mathbbm{X}=
\begin{pmatrix}
x_1\\
x_2\\
\vdots\\
x_n
\end{pmatrix},
\quad
\mathbbm{X}'=
\begin{pmatrix}
x_1'\\
x_2'\\
\vdots\\
x_n'
\end{pmatrix}.
\end{equation}
The two coordinates are related as $x_i'=\sum_{k=1}^n R_{ik}x_{k}$:
\begin{equation}
\mathbbm{X}'=\mathbbm{R}\mathbbm{X},
\label{eq:Xp=Rx}
\end{equation}
where the transformation matrix 
$\mathbbm{R}$ is an $n\times n$ invertible matrix with
known elements.
Cramer's rule states that the $i$th element of $\mathbbm{X}$ is determined as
\begin{equation}
\label{cramer-rule}
x_i=\frac{\mathscr{Det}[\mathbbm{R}^{(i)}(\mathbbm{X}')]}{\mathscr{Det}[\mathbbm{R}]},
\end{equation}
where 
$\mathscr{Det}$ represents the determinant,
$x_i$ is the $i$th element of $\mathbbm{X}$,
and $\mathbbm{R}^{(i)}(\mathbbm{X}')$ is the $n\times n$ matrix 
identical to $\mathbbm{R}$ except that
the $i$th column is replaced with $\mathbbm{X}'.$

We define the $j\times j$ square matrix $\mathbbm{R}_{[j\times j]}$ whose  elements are
identical to those of $\mathbbm{R}$ for $j=1$, $\cdots$, $n$:
\begin{equation}
\label{Rjj}
\mathbbm{R}_{[j\times j]}=
\begin{pmatrix}
R_{11}&\cdots &R_{1j}\\
\vdots&\ddots& \vdots\\
R_{j1}&\cdots &R_{jj}
\end{pmatrix}.
\end{equation}
Thus $\mathbbm{R}_{[n\times n]}=\mathbbm{R}$.
We define the column vectors $\mathbbm{X}'_{[j]}$ 
and $\mathbbm{R}^{(k)}_{[j]}$ with $j$ rows as
\begin{equation}
\label{XPj}
\mathbbm{X}'_{[j]}=
\begin{pmatrix}
x_1'\\\vdots\\x_j'
\end{pmatrix},
\quad
\mathbbm{R}^{(k)}_{[j]} =
\begin{pmatrix}
R_{1k}\\
\vdots\\
R_{jk}
\end{pmatrix}
.
\end{equation}
We define the partial contribution of $\mathbbm{X}'_{[j]}$
originated from the higher-dimensional coordinates $x_{j+1}$ through $x_n$ as
\begin{equation}
\label{XperpJ}
 \mathbbm{X}'_{\perp[j]}
 = \sum_{k=j+1}^{n} 
x_k \begin{pmatrix}
 R_{1k}\\
\vdots
\\
R_{jk}
\end{pmatrix}.
\end{equation}
For $j=n$, we set $\mathbbm{X}'_{\perp[n]}=\mathbb{0}$, the $n$-dimensional null vector.

We denote $\mathscr{D}_j$ by  the determinant of $\mathbbm{R}_{[j\times j]}$:
\begin{eqnarray}
\label{Detj}
\mathscr{D}_j
&=&\mathscr{Det}[\mathbbm{R}_{[j\times j]}]
\nonumber\\
&=&
\sum_{i_1=1}^{j}
\cdots
\sum_{i_j=1}^{j}
\epsilon_{i_1\cdots i_j}R_{1i_1}R_{2i_2}\cdots R_{ji_j},
\end{eqnarray}
where $\epsilon_{i_1\cdots i_j}$ is the $j$-dimensional
totally antisymmetric Levi-Civita symbol.

We denote $\mathbbm{R}^{(i)}_{[j\times j]}(\mathbbm{V})$ by the $j\times j$ square matrix 
identical to $\mathbbm{R}_{[j\times j]}$ defined in Eq.~(\ref{Rjj}) except that the $i$th column is replaced with
a column vector $\mathbbm{V}$ with $j$ elements:
\begin{equation}
\label{RJJV}
\mathbbm{R}^{(i)}_{[j\times j]}(\mathbbm{V})
=
\begin{pmatrix}
\mathbbm{R}^{(1)}_{[j]}\cdots
\mathbbm{R}^{(i-1)}_{[j]}\,
\mathbbm{V}\,
\mathbbm{R}^{(i+1)}_{[j]}\cdots
\mathbbm{R}^{(j)}_{[j]} 
\end{pmatrix}
,
\end{equation}
where $\mathbbm{R}^{(i)}_{[n\times n]}(\mathbbm{V})=\mathbbm{R}^{(i)}(\mathbbm{V})$.

\vspace{0.3cm}

\section{Proof of Cramer's Rule \label{sec;proof}}
In this section, we illustrate the proof of Cramer's rule with the aid of
the Dirac delta function.
As a byproduct, we derive a generalized version of
Cramer's rule involving a partial set of variables.

We multiply the unities
\begin{equation}
1=\int_{-\infty}^\infty dx'_i\,\delta(x_i'- \sum_{k=1}^n R_{ik}x_{k}) 
\end{equation}
to $\mathbbm{X}$ for $i=1$ through $n$.
Then we find that
\begin{equation}
\label{X}
\mathbbm{X}=
\int_{-\infty}^\infty dx'_1
 \cdots
\int_{-\infty}^\infty dx'_n
\mathbbm{X}\prod_{k=1}^n
\delta(x_k'-\sum_{a_k=1}^nR_{ka_k}x_{a_k}).
\end{equation}
The multiple integration over $x_1'$, $\cdots$, $x_n'$ can be replaced with 
the integration over the variables 
$x_1$, $\cdots$, $x_n$ as
\begin{equation}
\label{X-un}
\int_{-\infty}^\infty dx'_{1}
 \cdots
\int_{-\infty}^\infty dx'_n
=
\int_{-\infty}^\infty dx_1
 \cdots
\int_{-\infty}^\infty dx_n
\mathscr{J},
\end{equation}
where 
$\mathscr{J}$ is the Jacobian
\begin{equation}
\mathscr{J}
=
\left|
 \mathscr{Det}
 \begin{pmatrix}
\frac{\partial x_1'}{\partial x_1}&\cdots&\frac{\partial x_1'}{\partial x_n}\\
\vdots&\ddots&\vdots\\
\frac{\partial x_n'}{\partial x_1}&\cdots&\frac{\partial x_n'}{\partial x_n} 
 \end{pmatrix}
\right|
=
\left|
 \mathscr{Det}[\mathbbm{R}]
\right|
=
|\mathscr{D}_n|
.
\end{equation}

Then the integral (\ref{X}) reduces into
\begin{eqnarray}
\label{X-un-after}
\mathbbm{X} = 
\int_{-\infty}^\infty \!\!dx_1
 \cdots
\int_{-\infty}^\infty \!\!dx_n
|\mathscr{D}_n|
\mathbbm{X}
\prod_{k=1}^n
\delta(x_k'-\sum_{a_k=1}^nR_{ka_k}x_{a_k}).
\nonumber\\
\end{eqnarray}
The integration over $x_1$ in Eq.~(\ref{X-un-after}) can be carried out as
\begin{eqnarray}
\label{after-x1}
&&
\int_{-\infty}^\infty dx_1|\mathscr{D}_n|\mathbbm{X}\delta(x_1'-\sum_{ k =1}^nR_{1k}x_{k})
\nonumber\\
&=&
\int_{-\infty}^\infty dx_1
 |\mathscr{D}_n| 
 \mathbbm{X}
\delta\left[
\mathscr{Det}[
\mathbbm{R}^{(1)}_{[1\times 1]}(\mathbbm{X}'_{[1]}-\mathbbm{X}'_{\perp[1]})]
-\mathscr{D}_1 x_1 
\right]
\nonumber\\
&=&
\int_{-\infty}^\infty dx_1
 \frac{|\mathscr{D}_n|} {|\mathscr{D}_1|} 
 \mathbbm{X}
\delta\left[x_1-
\frac{\mathscr{Det}[
\mathbbm{R}^{(1)}_{[1\times 1]}(\mathbbm{X}'_{[1]}-\mathbbm{X}'_{\perp[1]})]}{\mathscr{D}_1}
\right]
\nonumber\\
&=&
\frac{|\mathscr{D}_n|}{|\mathscr{D}_1|}
\mathbbm{X}^{\{1\}},
\end{eqnarray}
where
$\mathscr{D}_1=\mathscr{Det}[\mathbbm{R}_{[1\times 1]}]=R_{11}$ and
\begin{equation}
\label{eq:X1-matrix}
\mathbbm{X}^{\{1\}}\equiv
\begin{pmatrix}
x_1(x_1',x_{2},\cdots,x_n)\\
x_2
\\
\vdots\\
x_n
\end{pmatrix}.
\end{equation}
The notation $\mathbbm{X}^{\{1\}}$ indicates that
$x_1$ in $\mathbbm{X}$ is replaced with the expression  in terms of
$x_1'$ and $x_k$ for $k=2$, $\cdots$, $n$ as
\begin{eqnarray}
\label{x1}
x_1&=&\frac{1}{\mathscr{D}_1}(x_1'- \sum_{k=2}^nR_{1k}x_k)
\nonumber\\
&=&\frac{1}{\mathscr{D}_1}
\mathscr{Det}[\mathbbm{R}^{(1)}_{[1\times 1]}(\mathbbm{X}'_{[1]}-\mathbbm{X}'_{\perp[1]})].
\end{eqnarray}
We shall find that, after carrying out the integrations over
$x_1$ through $x_j$, these coordinates are replaced with the corresponding
linear combinations of $x_1'$ through $x_j'$ and $x_{j+1}$ through $x_n$.
Thus we define $\mathbbm{X}^{\{j\}}$ by the generalized version of $\mathbbm{X}^{\{1\}}$
for $j=1,\cdots,n$ as
\begin{equation}
\label{XJ}
\mathbbm{X}^{\{j\}}=
\begin{pmatrix}
x_1(x_1',\cdots,x_j',x_{j+1},\cdots,x_n)\\
\vdots
\\
x_j(x_1',\cdots,x_j',x_{j+1},\cdots,x_n)
\\
x_{j+1}
\\
\vdots\\
x_n
\end{pmatrix}.
\end{equation}
\begin{widetext}
In a similar manner, the integration over $x_2$ can be carried out as
\begin{eqnarray}
\label{after-x2}
\int_{-\infty}^\infty dx_2
\frac{|\mathscr{D}_n|}{|\mathscr{D}_1|}
\mathbbm{X}^{\{1\}}
\delta(x_2'-\sum_{k=1}^nR_{2k}x_k) 
&=&
\int_{-\infty}^\infty dx_2
 \frac{|\mathscr{D}_n|}{|\mathscr{D}_2|} 
\mathbbm{X}^{\{1\}}
\delta\left[
\frac{1}{\mathscr{D}_2}
\left(
\mathscr{Det}[
\mathbbm{R}^{(2)}_{[2\times 2]}(\mathbbm{X}'_{[2]}-\mathbbm{X}'_{\perp[2]})]
\right)-x_2
\right]
\nonumber\\
&=&
\frac{|\mathscr{D}_n|}{|\mathscr{D}_2|}\mathbbm{X}^{\{2\}}.
\end{eqnarray}
Because the derivation of Eq.~\eqref{after-x2} is rather tricky,
we explain how we have obtained the results in detail:
\begin{eqnarray}
\label{detail-x2}
\left(
x_2'-\sum_{k=1}^nR_{2k}x_k\right)
&=&
\frac{1}{\mathscr{D}_1}
\left(
R_{11}x_2'
-R_{11}R_{21}x_1
-R_{11}\sum_{k=2}^nR_{2k}x_k 
\right)
\nonumber\\
&=&
\frac{1}{\mathscr{D}_1}
\left[
(R_{11}x_2'
- R_{21}x'_1)
-\sum_{k=3}^n
\left(
R_{11}R_{2k}-R_{21}R_{1k}
\right)x_k 
-
\left(
R_{11}R_{22}-R_{21}R_{12}
\right)x_2
\right]
\nonumber\\
&=&
\frac{1}{\mathscr{D}_1}
\left[
\mathscr{Det} 
\begin{pmatrix}
R_{11}&x_1'\\
R_{21}&x_2'
\end{pmatrix}
-\sum_{k=3}^n\mathscr{Det} 
\begin{pmatrix}
R_{11}&R_{1k}\\
R_{21}&R_{2k}
\end{pmatrix}x_k
-\mathscr{Det} 
\begin{pmatrix}
R_{11}&R_{12}\\
R_{21}&R_{22}
\end{pmatrix}x_2
\right]
\nonumber\\
&=&
\frac{1}{\mathscr{D}_1}
\left[
\mathscr{Det}[
\mathbbm{R}^{(2)}_{[2\times 2]}(\mathbbm{X}'_{[2]})]
-
\mathscr{Det}[
\mathbbm{R}^{(2)}_{[2\times 2]}(\mathbbm{X}'_{\perp[2]})]-\mathscr{D}_2 x_2
\right]
\nonumber\\
&=&
\frac{\mathscr{D}_2}{\mathscr{D}_1}
\left[
\frac{1}{\mathscr{D}_2}
\mathscr{Det}[
\mathbbm{R}^{(2)}_{[2\times 2]}(\mathbbm{X}'_{[2]}-\mathbbm{X}'_{\perp[2]})]
- x_2
\right].
\end{eqnarray}
On the right side of the first line of Eq.~\eqref{detail-x2},
we have imposed $\mathscr{D}_1=R_{11}$ from Eq.~\eqref{Detj} and
pulled out the $k=1$ term from the summation over $k$.
According to Eqs.~\eqref{eq:X1-matrix} and \eqref{x1},
we replace $x_1$ with a linear combination of $x_1'$ and $x_2,\cdots, x_n$.
Collecting the primed-coordinate, $x_2$, and unprimed-coordinate
contributions of $k\ge 3$ separately, we obtain the second line of Eq.~\eqref{detail-x2}.
Then the remaining derivation 
from the third to the last line of Eq.~\eqref{detail-x2} is straightforward.

After evaluating the $x_2$ integral, we end up with
$\mathbbm{X}^{\{2\}}$
defined in  Eq.~\eqref{XJ} up to an overall factor of determinant ratio $|\mathscr{D}_n|/|\mathscr{D}_2|$.
After the double integrations over $x_1$ in Eq.~(\ref{after-x1}) and over $x_2$ in Eq.~(\ref{after-x2}), 
$x_1$ and $x_2$ are expressed in terms of $x_1'$, $x_2'$, and $x_k$'s for $k=3$, $\cdots$, $n$ as
\begin{equation}
x_i=\frac{1}{\mathscr{D}_2}
\mathscr{Det}[
\mathbbm{R}^{(i)}_{[2\times 2]}(\mathbbm{X}'_{[2]}-\mathbbm{X}'_{\perp[2]})],\quad i=1,\,2,
\end{equation}
where the expression for $x_1$ is obtained
by inserting $x_2$ into Eq.~(\ref{x1}).

In Appendix~\ref{app:mathematical-induction}, we show that
the integration over $x_j$ can be carried out
by employing mathematical induction as
\begin{eqnarray}
\label{after-xj}
\int_{-\infty}^\infty\!\!\!dx_j
\frac{|\mathscr{D}_n|}{|\mathscr{D}_{j-1}|}
\mathbbm{X}^{\{j-1\}}
\delta(x_j'-\sum_{k=1}^nR_{jk}x_k) 
&=&
\int_{-\infty}^\infty\!\!\! dx_j
\frac{|\mathscr{D}_n|}{|\mathscr{D}_{j-1}|}
\mathbbm{X}^{\{j-1\}}
\delta\left[
\frac{1}{\mathscr{D}_{j-1}}\left(
\mathscr{Det}[
\mathbbm{R}^{(j)}_{[j\times j]}(\mathbbm{X}'_{[j]}-\mathbbm{X}'_{\perp[j]})]
-\mathscr{D}_j x_j \right)
\right]
\nonumber\\
&=&
\int_{-\infty}^\infty\!\!\! dx_j
\frac{|\mathscr{D}_n|}{|\mathscr{D}_{j}|}
\mathbbm{X}^{\{j-1\}}
\delta\left[x_j-
\frac{\mathscr{Det}[
\mathbbm{R}^{(j)}_{[j\times j]}(\mathbbm{X}'_{[j]}-\mathbbm{X}'_{\perp[j]})]}{\mathscr{D}_j}
\right]
\nonumber\\
&=&
\frac{|\mathscr{D}_n|}{|\mathscr{D}_{j}|}
\mathbbm{X}^{\{j\}},
\label{intj}%
\end{eqnarray}
where $\mathbbm{R}^{(i)}_{[j\times j]}(\mathbbm{V})$ and
$\mathbbm{X}'_{\perp[j]}$ are defined in Eqs.~(\ref{RJJV}) and (\ref{XperpJ}),
respectively. We note that $\mathbbm{X}'_{\perp[j]}$ depends only on
$x_{j+1}$ through $x_n$.

After the multiple integrations over $x_1$ in Eq.~(\ref{after-x1}) through over $x_j$ in Eq.~(\ref{after-xj}), 
$x_i$ for $i=1$ through $j$ are expressed in terms of $x_1'$, $\cdots$, $x_j'$
and $x_k$'s for $k=j+1$, $\cdots$, $n$ as
\begin{equation}
\label{xij}
x_i
=\frac{1}{\mathscr{D}_j}
\mathscr{Det}[\mathbbm{R}^{(i)}_{[j\times j]}(\mathbbm{X}'_{[j]}-\mathbbm{X}'_{\perp[j]})], \quad i=1,\,2,\cdots,\, j.
\end{equation}
Equation (\ref{xij}) implies that a similar form to Cramer's rule holds in the change of a partial set of variables.
It turns out that Cramer's rule is actually a special case for $j=n$
of the general formula~(\ref{xij}),
which we call \textit{Cramer's rule for a partial set of variables
or coordinates}.
To our best knowledge, the generalized formula (\ref{xij}) is new.
In the following section, we demonstrate how the rule (\ref{xij})
applies with a simple mechanical system.

We repeat the same procedure until we reach $j=n$ at which $\mathbbm{X}_{\perp[n]}$ vanishes:
\begin{eqnarray}
\int_{-\infty}^\infty dx_n\frac{|\mathscr{D}_n|}{|\mathscr{D}_{n-1}|}\mathbbm{X}^{\{n-1\}}
\delta(x_n'- \sum_{k=1}^nR_{nk}x_k) 
&=&
\int_{-\infty}^\infty\!\!\! dx_n
\frac{|\mathscr{D}_n|}{|\mathscr{D}_{n-1}|}
\mathbbm{X}^{\{n-1\}}
\delta\left[
\frac{1}{\mathscr{D}_{n-1}}\left(
\mathscr{Det}[
\mathbbm{R}^{(n)}(\mathbbm{X}')]
- 
\mathscr{D}_n
x_n\right)
\right]
\nonumber\\
&=&
\int_{-\infty}^\infty\!\!\! dx_n
\mathbbm{X}^{\{n-1\}}
\delta\left[x_n-\frac{\mathscr{Det}[\mathbbm{R}^{(n)}(\mathbbm{X}')]}{\mathscr{Det}[\mathbbm{R}]}
\right]
\nonumber\\
&=&
\mathbbm{X}^{\{n\}}.
\end{eqnarray}
As a result, for $i=1$, $\cdots$, $n$,
\begin{eqnarray}
\label{xi}
x_i= \frac{\mathscr{Det}[\mathbbm{R}^{(i)}(\mathbbm{X}')]}{\mathscr{Det}[\mathbbm{R}]},
\end{eqnarray}
where $\mathbbm{R}^{(i)}(\mathbbm{X}')=\mathbbm{R}^{(i)}_{[n\times n]}(\mathbbm{X}'_{[n]})$.
This completes the proof of Cramer's rule~(\ref{cramer-rule}).

\section{Application \label{sec;app}}
As an example of the physical application of 
Cramer's rule for a partial set of variables
in Eq.~(\ref{xij}),
we consider the accelerating motion of 
a series of $n$ blocks connected with ropes 
of negligible mass as shown in Fig.~\ref{figure:boxes}.
\begin{figure}
\centering
\includegraphics[width=0.8\columnwidth]{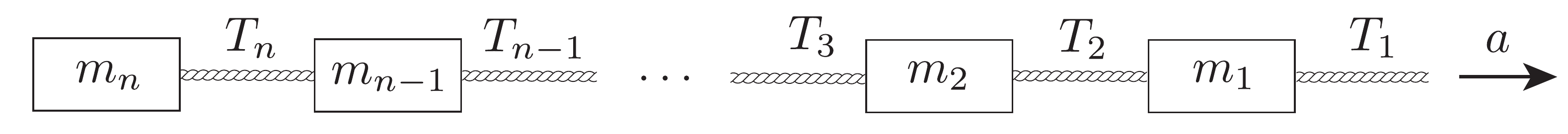}
\caption{\label{figure:boxes}
A series of $n$ blocks connected with ropes 
of negligible mass moving with the acceleration $a$.
$T_i$ is the tension of the rope connecting the blocks
of mass $m_i$ and $m_{i-1}$.}
\end{figure}
We assume that every rope is tight and do not consider the restoring 
force of the rope.
Newton's second law states that
\begin{eqnarray}
\label{eq:tension-force-example}
a
\begin{pmatrix}
m_1
\\
m_2
\\
m_3
\\
\vdots
\\
m_{n-2}
\\
m_{n-1}
\\
m_n
\end{pmatrix}
=
\begin{pmatrix*}[r]
1 & -1 & 0  & \cdots & 0 & 0 & 0 
\\
0 & 1 & -1  & \cdots & 0 & 0 & 0
\\
0 & 0 & 1  & \cdots & 0 & 0 & 0
\\
\vdots & \vdots  &\vdots &\ddots &\vdots & \vdots & \vdots
\\
0 & 0 &  0 &\cdots & 1 & -1 & 0
\\
0 & 0 &  0 &\cdots & 0 & 1 & -1
\\
0 & 0 &  0 &\cdots & 0 & 0 & 1
\end{pmatrix*}
\begin{pmatrix}
T_1
\\
T_2
\\
T_3
\\
\vdots
\\
T_{n-2}
\\
T_{n-1}
\\
T_{n}
\end{pmatrix},
\end{eqnarray}
where $a$ is the acceleration of the system,
$m_i$ is the mass of the $i$th block, and $T_i$ is the tension 
of the rope between $m_i$ and $m_{i-1}$.
$\mathbbm{X}'$, $\mathbbm{R}$, and $\mathbbm{X}$ in Eq.~(\ref{eq:Xp=Rx})
can be read off from Eq.~(\ref{eq:tension-force-example})
to find that
\begin{equation}
x_i = T_i,
\quad\textrm{and}\quad
x_i' =  m_i a,
\quad
\textrm{for $i=1,\cdots,n$.}
\end{equation}

Let us choose the partial set of variables
as $T_1$, $T_2$, and $T_3$ with $j=3$.
By making use of Cramer's rule for a partial set of variables
in Eq.~(\ref{xij}), we can express
$T_1$, $T_2$, and $T_3$
in terms of $m_1a$, $m_2a$, and $m_3a$ and $T_4$, $\cdots$, $T_n$ as follows: 
\begin{eqnarray}
\label{eq:Ti-cramer}
T_1&=&
\frac{1}{\mathscr{D}_3}
\mathscr{Det}[\mathbbm{R}^{(1)}_{[3\times 3]}(\mathbbm{X}'_{[3]}-\mathbbm{X}'_{\perp[3]})],
\nonumber \\
T_2&=&
\frac{1}{\mathscr{D}_3}
\mathscr{Det}[\mathbbm{R}^{(2)}_{[3\times 3]}(\mathbbm{X}'_{[3]}-\mathbbm{X}'_{\perp[3]})],
\nonumber \\
T_3&=&
\frac{1}{\mathscr{D}_3}
\mathscr{Det}[\mathbbm{R}^{(3)}_{[3\times 3]}(\mathbbm{X}'_{[3]}-\mathbbm{X}'_{\perp[3]})],
\end{eqnarray}
where 
\begin{subequations}
\label{eq:Ti-derivation}
\begin{eqnarray}
\mathscr{D}_3
&=&
\mathscr{Det}[\mathbbm{R}_{[3\times 3]}]
=
\mathscr{Det}
\begin{pmatrix*}[r]
1 & -1 & 0   
\\
0 & 1 & -1 
\\
0 & 0 & 1 
\end{pmatrix*}
=1,
\nonumber \\
\mathbbm{X}'_{[3]}-\mathbbm{X}'_{\perp[3]}
&=&
\begin{pmatrix}
x_1'
\\
x_2'
\\
x_3'
\end{pmatrix}
-
\sum_{k=4}^n
x_k
\begin{pmatrix}
R_{1k}
\\
R_{2k}
\\
R_{3k}
\end{pmatrix}
=
\begin{pmatrix}
x_1'
\\
x_2'
\\
x_3'
\end{pmatrix}
-
x_4
\begin{pmatrix}
0
\\
0
\\
-1
\end{pmatrix}
=
\begin{pmatrix}
m_1a
\\
m_2a
\\
m_3a+T_4
\end{pmatrix},
\end{eqnarray}
and
\begin{eqnarray}
\mathscr{Det}[\mathbbm{R}^{(1)}_{[3\times 3]}(\mathbbm{X}'_{[3]}-\mathbbm{X}'_{\perp[3]})]
&=&
\mathscr{Det}
\begin{pmatrix}
m_1a & -1 & \phantom{-}0   
\\
m_2a & \phantom{-}1 & -1 
\\
m_3a+T_4& \phantom{-}0 & \phantom{-}1 
\end{pmatrix}
=
(m_1+m_2+m_3)a+T_4,
\nonumber \\
\mathscr{Det}[\mathbbm{R}^{(2)}_{[3\times 3]}(\mathbbm{X}'_{[3]}-\mathbbm{X}'_{\perp[3]})]
&=&
\mathscr{Det}
\begin{pmatrix}
1 & m_1a & \phantom{-}0   
\\
0 & m_2a & -1 
\\
0& m_3 a+T_4 & \phantom{-}1 
\end{pmatrix}
=
(m_2+m_3)a+T_4,
\nonumber \\
\mathscr{Det}[\mathbbm{R}^{(3)}_{[3\times 3]}(\mathbbm{X}'_{[3]}-\mathbbm{X}'_{\perp[3]})]
&=&
\mathscr{Det}
\begin{pmatrix}
1 & -1 & m_1a
\\
0 & \phantom{-}1 & m_2 a
\\
0& \phantom{-}0 & m_3 a+T_4
\end{pmatrix}
=
m_3a+T_4.
\end{eqnarray}
\end{subequations}
Substituting the values in Eqs.~(\ref{eq:Ti-derivation}) 
into Eq.~(\ref{eq:Ti-cramer}),
we obtain
\begin{eqnarray}
\label{eq:app-exam-1-result}
T_1&=&(m_1+m_2+m_3)a+T_4,
\nonumber \\
T_2&=&(m_2+m_3)a+T_4,
\nonumber \\
T_3&=&m_3a+T_4.
\end{eqnarray}
Therefore, using Eq.~(\ref{xij}), 
we could express the selected set of variables $T_1$, $T_2$, $T_3$
in terms of the transformed variables, $m_1a$, $m_2a$, $m_3a$, 
and could keep the remaining set of variables $T_k$ for $k\ge 4$
intact.
This partial transformation enabled by 
Cramer's rule for a partial set of variables
in Eq.~(\ref{xij}) could also be used in more general cases
involving a large number of variables some of which do not have to 
be specified
in terms of the transformed variables.

Recursive applications of the Cramer's rule
for a partial set of variables
in Eq.~(\ref{xij})
through $j=n$ lead to
\begin{eqnarray}
\begin{pmatrix}
T_1
\\
T_2
\\
T_3
\\
\vdots
\\
T_{n-2}
\\
T_{n-1}
\\
T_{n}
\end{pmatrix}
=
a
\begin{pmatrix*}[r]
1 & \phantom{-}1 & \phantom{-}1  & \cdots & 1 & \phantom{-}1 & \phantom{-} 1
\\
0 & 1 & 1  & \cdots & 1 & 1 & 1
\\
0 & 0 & 1  & \cdots & 1 & 1 & 1
\\
\vdots & \vdots & \vdots &\ddots &\vdots &\vdots & \vdots 
\\
0 & 0 &  0 &\cdots & 1 & 1 & 1
\\
0 & 0 &  0 &\cdots & 0 & 1 & 1
\\
0 & 0 &  0 &\cdots & 0 & 0 & 1
\end{pmatrix*}
\begin{pmatrix}
m_1
\\
m_2
\\
m_3
\\
\vdots
\\
m_{n-2}
\\
m_{n-1}
\\
m_n
\end{pmatrix}.
\end{eqnarray}
As a result, $T_i$'s are completely solved as 
\begin{equation}
\label{eq:Ti-complete}
T_i = \left(\sum_{k=i}^n m_i\right)a.
\end{equation}
While the solution (\ref{eq:Ti-complete}) provides complete
information of every tension, the solution (\ref{eq:app-exam-1-result})
for a partial set of variables
in which one is interested 
contains every specific information that one wants to know. 
Thus the solution (\ref{eq:app-exam-1-result}) should be particularly  
economical as the number of variables becomes very large. 

\end{widetext}

\section{Discussion \label{sec;dis}}
We have derived Cramer's rule for a solution of the system of linear equations, $\mathbbm{X}^\prime = \mathbbm{R} \mathbbm{X}$,
where $\mathbbm{X}$ and $\mathbbm{X}^\prime$ are interpreted as $n$-dimensional Cartesian-coordinate column vectors
in Eq.~(\ref{xxp})
and $\mathbbm{R}$ is the corresponding transformation matrix.
The introduction of a trivial identity (\ref{X})  involving Dirac delta functions has leaded to a way to reach
the systematic transformation of the original coordinates
$x_i$'s
into new coordinates $x^\prime_i$'s in a straightforward manner resulting in 
the completion of the proof of Cramer's rule~(\ref{cramer-rule}).
To our best knowledge, this proof is new and, furthermore, could be useful to the pedagogical training of
recursive application of integration of Dirac delta functions. 

In an intermediate step, we have derived Eq.~(\ref{xij}), where $x_i$ ($i=1,\cdots,j $ ) is
expressed in terms of $x^\prime_i$ and $x_k$ ($k=j+1,\cdots,n$) with the transformation matrix $\mathbbm{R}$.
In order to find the solution, we have to know information for the elements of the submatrix $\mathbbm{R}_{[j\times j]}$
as well as those of another $j\times (n-j)$ matrix, $R_{ik}$.
We call Eq.~(\ref{xij}) as 
\textit{Cramer's rule for a partial set of variables or coordinates}.
This rule can be immediately applied to change a partial set of generalized coordinates of a mechanical system.
The original Cramer's rule in Eq.~(\ref{cramer-rule}) 
is actually a special case for $j=n$
of the general formula~(\ref{xij}).
This rule can
further be extended to any linear problems,
such as Kirchhoff's circuit laws, and any oscillating systems whose
displacements are small.

\begin{acknowledgments}
As members of the Korea Pragmatist Organization for Physics Education
(\textsl{KPOP}$\mathscr{E}$), 
the authors thank the remaining members of \textsl{KPOP}$\mathscr{E}$ 
for useful discussions.
This work is supported in part by the National Research Foundation
of Korea (NRF) under the BK21+ program at Korea University, 
\textit{Initiative for Creative and Independent Scientists},
and by grants funded by the Korea government (MSIT),
Grants No. NRF-2017R1A2B4011946 (C.Y.), No. NRF-2017R1E1A1A01074699 (J.L.)
and No. NRF-2020R1A2C3009918 (J.E.).
The work of C.Y. is also supported in part by Basic Science Research Program through the National Research Foundation of Korea(NRF) funded by the Ministry of Education(2020R1I1A1A01073770).

\end{acknowledgments}

\appendix

\begin{widetext}
\section{Proof of Eq.~(\ref{after-xj}) by
mathematical induction}
\label{app:mathematical-induction}
	
In this Appendix, we present a proof of Eqs.~\eqref{intj} and \eqref{xij}
by mathematical induction.
We have proved the statement for $j=1$ in Eqs.~(\ref{after-x1}) and 
(\ref{x1}).
We assume that Eqs.~\eqref{intj} and \eqref{xij}
are satisfied for a given $j=p-1$.
In order to prove that the counterparts of 
Eqs.~(\ref{after-x1}) and 
(\ref{x1}) for $j=p$ are true,
it is crucial to verify that the following relation, which involves 
the argument of the delta function in the left side of Eq.~\eqref{intj}, 
is true:
\begin{equation}
\mathscr{D}_{p-1} \left[
x_{p}'-\sum_{k=1}^nR_{pk}x_k
\right]
=
\mathscr{Det}[
\mathbbm{R}^{(p)}_{[p\times p]}(\mathbbm{X}'_{[p]}-\mathbbm{X}'_{\perp[p]})]
-\mathscr{D}_{p} x_{p}.
\label{a1}%
\end{equation}
The left side of Eq.~\eqref{a1} can be rewritten as
\begin{eqnarray}
\mathscr{D}_{p-1}\left[
x_{p}'-\sum_{k=1}^nR_{p k}x_k
\right]
&=&
\mathscr{D}_{p-1}\left[
x_{p}'-\sum_{k=1}^{p-1} R_{p k}x_k
-R_{pp} x_p
-\sum_{k=p+1}^{n} R_{p k}x_k
\right]
\nonumber \\
&=&
\mathscr{D}_{p-1}\left[
x_{p}'-\sum_{k=1}^{p-1} R_{p k}\left(
\frac{1}{\mathscr{D}_{p-1}}
\mathscr{Det}[\mathbbm{R}^{(k)}_{[(p-1)\times (p-1)]}(\mathbbm{X}'_{[p-1]}-\mathbbm{X}'_{\perp[p-1]})]\right)
-R_{pp} x_p
-\sum_{k=p+1}^{n} R_{p k}x_k
\right]
\nonumber \\
&=&
A(x_1',\cdots,x_p')+B(x_p)+C(x_{p+1},\cdots,x_n),
\label{a2}
\end{eqnarray}
where $x_k$ for $k=1,\cdots,p-1$ in the second line is replaced with
the relation \eqref{xij}.
In the last line, we collect the primed-coordinate, $x_p$, and
unprimed-coordinate contributions of $k\ge p+1$ separately,
which are defined by functions, $A(x_1',\cdots,x_p')$, $B(x_p)$, and $C(x_{p+1},\cdots,x_n)$, respectively.

If we list up the matrix representation for 
$\mathbbm{R}^{(k)}_{[(p-1)\times (p-1)]}
(\mathbbm{X}'_{[p-1]}-\mathbbm{X}'_{\perp[p-1]})$, then
we find that
\begin{equation}
\label{DET-MAT}
\mathbbm{R}^{(k)}_{[(p-1)\times (p-1)]}(\mathbbm{X}'_{[p-1]}-\mathbbm{X}'_{\perp[p-1]}) 
=
\begin{pmatrix}
R_{11} & \cdots & R_{1\,k-1} & x_1' - \sum_{\ell=p}^{n} R_{1\ell} x_{\ell} & \cdots & R_{1\, p-1} \\
R_{21} & \cdots &  R_{2\,k-1} &x_2'- \sum_{\ell=p}^{n} R_{2\ell} x_{\ell}&  \cdots & R_{2\, p-1} \\
\vdots & \ddots& \vdots & \vdots & \ddots & \vdots \\
R_{p-1\,1} & \cdots  &R_{p-1\,k-1} &x_{p-1}'- \sum_{\ell=p}^{n} R_{p-1\,\ell} x_{\ell}& \cdots & R_{p-1\,p-1}
\end{pmatrix}.\phantom{X}
\end{equation}
The primed-coordinate contribution $A(x_1',\cdots,x_p')$ can be collected
from Eqs.~\eqref{a2} and \eqref{DET-MAT} as
\begin{eqnarray}
\label{A}
A(x_1',\cdots,x_p')
&=& 
\mathscr{D}_{p-1}x_{p}'
-
\sum_{k=1}^{p-1} R_{p k}\,\mathscr{Det}
\begin{pmatrix}
R_{11} & \cdots & R_{1\,k-1} & x_1'  & R_{1\,k+1} & \cdots & R_{1\,p-1} \\
R_{21} & \cdots &  R_{2\,k-1} &x_2'& R_{2\,k+1}&  \cdots & R_{2\,p-1} \\
\vdots & \ddots& \vdots & \vdots &\vdots& \ddots & \vdots \\
R_{p-1\,1} & \cdots  &R_{p-1\,k-1} &x_{p-1}'& R_{p-1\,k+1} & \cdots & R_{p-1\, p-1}
\end{pmatrix}
\nonumber\\
&=& 
\mathscr{Det}
\begin{pmatrix}
R_{11} & R_{12} & \cdots & R_{1\,p-1} & x_1' \\
R_{21} & R_{22} & \cdots & R_{2\,p-1} & x_2' \\
\vdots & \vdots & \ddots & \vdots     & \vdots \\
R_{p1} & R_{p2} & \cdots & R_{p\,p-1} & x_p'
\end{pmatrix}
\nonumber \\
&=&
\mathscr{Det}[\mathbbm{R}^{(p)}_{[p\times p]}(\mathbbm{X}'_{[p]})],
\end{eqnarray}
where we have made use of  the following two properties of the determinant:  (i)
$\mathscr{Det}[\mathbbm{R}^{(k)}_{[p\times p]}(a_1\mathbbm{V}_1+a_2\mathbbm{V}_2 )]
=a_1\mathscr{Det}[\mathbbm{R}^{(k)}_{[p\times p]}(\mathbbm{V}_1)]
+a_2\mathscr{Det}[\mathbbm{R}^{(k)}_{[p\times p]}(\mathbbm{V}_2)]$
for any numbers $a_1$ and $a_2$.
(ii) An exchange of two distinct columns of a matrix flips the sign of the determinant.

In a similar manner, $B(x_p)$ and $C(x_{p+1},\cdots,x_n)$ can be collected
from Eqs.~\eqref{a2} and \eqref{DET-MAT} and simplified as
\begin{eqnarray}
B(x_p)
&=&
\sum_{k=1}^{p-1} R_{p k}\mathscr{Det}
\begin{pmatrix}
R_{11} & \cdots & R_{1\,k-1} &   R_{1p} x_p & R_{1\,k+1} & \cdots & R_{1\,p-1} \\
R_{21} & \cdots &  R_{2\,k-1} & R_{2p} x_p&   R_{2\,k+1} &\cdots & R_{2\,p-1} \\
\vdots & \ddots& \vdots & \vdots & \vdots & \ddots &\vdots \\
R_{p-1\,1} & \cdots  &R_{p-1\,k-1} &  R_{p-1\,p} x_p&  R_{p-1\,k+1} &\cdots & R_{p-1\, p-1}
\end{pmatrix}
-R_{pp}\mathscr{D}_{p-1} x_p
\nonumber \\
&=&
-x_p \mathscr{Det}
\begin{pmatrix}
R_{11} & R_{12} & \cdots & R_{1\,p-1} & R_{1p} \\
R_{21} & R_{22} & \cdots & R_{2\,p-1} & R_{2p} \\
\vdots & \vdots & \ddots & \vdots     & \vdots \\
R_{p1} & R_{p2} & \cdots & R_{p\,p-1} & R_{pp}
\end{pmatrix}
\nonumber \\
&=&
-\mathscr{D}_p x_p,
\label{B}
\\
C(x_{p+1},\cdots,x_n)
&=&
\sum_{k=1}^{p-1} R_{p k}\mathscr{Det}
\begin{pmatrix}
R_{11} & \cdots & R_{1\,k-1} &   \sum_{\ell=p+1}^{n} R_{1\ell} x_\ell & R_{1\,k+1} & \cdots & R_{1\,p-1} \\
R_{21} & \cdots &  R_{2\,k-1} & \sum_{\ell=p+1}^{n} R_{2\ell} x_\ell&   R_{2\,k+1} &\cdots & R_{2\,p-1} \\
\vdots & \ddots& \vdots & \vdots & \vdots & \ddots &\vdots \\
R_{p-1\,1} & \cdots  &R_{p-1\,k-1} &  \sum_{\ell=p+1}^{n} R_{p-1\,\ell} x_\ell&  R_{p-1\,k+1} &\cdots & R_{p-1\, p-1}
\end{pmatrix}
\nonumber\\&&
-
\mathscr{D}_{p-1}\left(
\sum_{\ell=p+1}^{n} R_{p \ell}x_\ell
\right)
\nonumber \\
&=&
- \mathscr{Det}
\begin{pmatrix}
R_{11} & R_{12} & \cdots & R_{1\,p-1} & \sum_{\ell=p+1}^{n} R_{1 \ell}x_\ell \\
R_{21} & R_{22} & \cdots & R_{2\,p-1} &\sum_{\ell=p+1}^{n} R_{2 \ell}x_\ell \\
\vdots & \vdots & \ddots & \vdots     & \vdots \\
R_{p1} & R_{p2} & \cdots & R_{p\,p-1} & \sum_{\ell=p+1}^{n} R_{p \ell}x_\ell
\end{pmatrix}
\nonumber \\
&=&
-\mathscr{Det}[
\mathbbm{R}^{(p)}_{[p\times p]}(\mathbbm{X}'_{\perp[p]})].
\label{C}%
\end{eqnarray}
Substituting Eqs.~\eqref{A}, \eqref{B}, and \eqref{C} into  \eqref{a2},
we obtain \eqref{a1}. 
As a result, the relations  \eqref{intj} and \eqref{xij} hold for any $j$ for $j=1$
through $n$ by mathematical induction.

\end{widetext}

\end{document}